\documentclass[preprint,3p]{elsarticle}

\usepackage{amsmath,amssymb,amsfonts,bbm,amsthm}
\usepackage{color}

\newtheorem{thm}{Theorem}

\newtheorem{lmm}{Lemma}[section]

\newtheorem{corr}{Corollary}

\newcommand{\Pro}{\mathbf{P}}
\newcommand{\Es}{{\mbox{\bf E}}}
\newcommand{\R}{\mathbbm{R}}

\newcommand{\K}{\mathbbm{K}}

\newcommand{\Hu}{H}    													

\newcommand{\cL}{\mathcal{L}}

\newcommand{\MLE}{\hat{\vartheta}}

\newtheorem{rem}{Remark}

\newcommand{\Rg}       {{\hbox{I\kern-.22em\hbox{R}}}}
\newcommand{\Pg}       {{\hbox{I\kern-.22em\hbox{P}}}}
\newcommand{\Eg}       {{\hbox{I\kern-.22em\hbox{E}}}}

\newcommand{\tr}   {\mbox{ \rm{tr} }}

\journal{Physica A}

\begin{document}

\begin{frontmatter}



\title{Experiment design for controlled partially observed fractional diffusion process}


\author{Chunhao Cai and Wujun Lv}

\address{Shanghai University of Economics and Finance, Shanghai, China}

\begin{abstract}
We consider a controlled second order differential equation which is partially observed with an additional fractional noise. we study the asymptotic (for large observation time) design problem of the input and give an efficient estimator of the unknown signal drift parameter. When the input depends on the unknow parameter, we will try the one-step estimation procedure using the Newton-Raphon method.

\end{abstract}


\end{frontmatter}


\section{Introduction}

\subsection{Historical survey}

Over the last decades the experiment design has been given a great deal of interest from the early statistics literature (see \textit{e.g.}
 \cite{kiefer, wald, whittele}) as well as in the engineering literature (see \textit{e.g.} \cite{gevers,goodwin3,goodwin}).

The experiment design consists two problem or two procedure: the first is to find the energy constraint of the input which can maximize the Fisher information. The second problem is under this input how to find an adaptive estimator. In this area, there are several approaches like sequential design and Bayesian design (see \textit{e.g.} \cite{goodwin, Lev, Meh} and the references therein).

We will  also find some works which concern on the partially observed models such as \cite{Aoki, Levin, Meh, Meh2, Oss}, where  linear signal - observation model perturbed by the white noise has been considered.

On the other hand, large sample asymptotic properties (the consistency and the asymptotical normality) of the Maximum Likelihood Estimator (MLE) with the fractional noise  \cite{cheridito,Klep02a,Coutin,BK10,Lototsky} have been got enough attention.

Some models of the experiment design with the fractional noise have been studied by Brouste, Cai, Kleptsyna and Popier \cite{BC13,BKP10b,BKP10}. In these works the optimal input what we have found does not depend the unknow parameters, that is to say it is very easy to obtain directly the Maximum Likelihood Estimator. In this paper, even some technical methods will be the same of the pevious works, we will consider the situation of complex-valued equation and in this case we will meet a very different problem in the estimation procedure: the optimal input will depend on the parameter. In this sense, we will use one-step procedure of estimation using the Newton-Raphson method.

The paper falls into four parts. In this introduction, we state our models and then we will give our main results in the second part. In the third part,we will try to do some transformation of the models and present the Newton-Raphson method. The proofs of two lemmas will end all our works.

\subsection{The Model And Statement Of The Problem}
We consider complex-valued functions $x(t)$, $u(t)$,$t\geq 0$ and a process $Y=(Y_t,t\geq 0)$, representing the signal and the observation respectively, governed by the following homogeneous linear system of ordinary and stochastic differential equations interpreted as integral equations:
\begin{equation}\label{eq:sec model}
\left\{\begin{array}{ccrr}
 \frac{d^2x}{dt^2}+k\frac{dx}{dt}+\vartheta x&=&u(t) \,,&x(0)=0\,,\\[.2cm]
dY_t &=& x(t)dt+dV_t^H\,,&Y_0=0.
\end{array}\right.
\end{equation}
Here, $V^H=(V_t^H,\,t\geq 0)$ is normalized fBm with Hurst Index $H\in [\frac{1}{2},1)$ and the coefficient $\vartheta$ and $k$ are positive constants. System \eqref{eq:sec model} has a uniquely defined solution process $(x, Y)$ where $Y$ is Gaussian but neither Markovian nor a semi- martingale for  $H\neq \frac{1}{2}$.

Suppose that the parameter $\vartheta$ is unknown and is to be estimated given the observed trajectory $Y^T=(Y_t,\,0\leq t\leq T)$. For a fixed value of the parameter $\vartheta$, let $\mathbf{P}^{T}_\vartheta$ denote the probability
measure, induced by $(X^T,Y^T)$ on the function space ${\cal C}_{[0,T]}\times{\cal C}_{[0,T]}$ and
let ${\cal F}_t^Y$ be the natural filtration of $Y$, ${\cal F}_t^Y=\sigma\left( Y_s , 0 \leq s \leq t \right)$.

Let ${\cal L}(\vartheta,Y^{T})$ be the likelihood, {\it i.e.} the Radon-Nikodym derivative of
$\mathbf{P}^{T}_\vartheta$, restricted to ${\cal F}_T^Y$ with respect to some reference measure
on ${\cal C}_{[0,T]}$.  In this setting, Fisher information stands for :
$$
{\cal I}_T\left( \vartheta , u\right) = - \Es_\vartheta \frac{\partial^2}{\partial \vartheta^2} \ln \cL_T(\vartheta,Y^{T})\,.
$$

Let us denote $\mathcal{U}_T$ some functional space of controls, that is defined by equation \eqref{sec control_spaces} and \eqref{sec vtou}. Let us therefore note
$$
\mathcal{J}_T(\vartheta)=\sup_{u \in \mathcal{U}_T} \mathcal{I}_T(\vartheta,u).
$$
Our main goal is to find estimator $\overline{\vartheta}_T$ of the parameter $\vartheta$ which are asymptotically efficient in the sense that, for any compact $\K \subset \R^+_*$,
\begin{equation} \label{asymp}
\sup_{\vartheta \in \K} {\cal J}_T(\vartheta) \Es_\vartheta \left( \overline{\vartheta}_T - \vartheta \right)^2 =1+o(1) \,,
\end{equation}
as $T \rightarrow \infty$.

\section{Main Result}
In this section, we will divide two different cases,  we will get the optimal input and study the properties of the MLE.

\subsection{Case of $k^2\geq 2\vartheta$}
In this subsection we will consider only the case that $k^2\geq 2\vartheta$
\begin{thm}\label{sec opt case 1}
The asymptotical optimal input in the class of controls is $\mathcal{U}_T$ is $u_{opt}^1(t)=\frac{\kappa_H}{\sqrt{2\lambda}}t^{H-\frac{1}{2}}$ where
\begin{equation} \label{sec def_constants}
\kappa_H=2\Hu\Gamma\left(\frac{3}{2}-\Hu\right)\Gamma\left(\frac{1}{2}+\Hu\right)
\quad \mbox{and} \quad \lambda= \frac{\Hu \Gamma(3-2\Hu) \Gamma(\Hu+\frac{1}{2})}{2(1-\Hu)\Gamma(\frac{3}{2}-\Hu)},
\end{equation}
and $\Gamma$ stands for the Gamma function. Moreover,
$$\lim_{T \to + \infty} \frac{{\cal J}_T^1(\vartheta) }{T} =\mathcal{I}^1(\vartheta),$$
where
\begin{equation}\label{eq:sec fisher 1}
\mathcal{I}^1(\vartheta)=\frac{1}{\vartheta^4}.
\end{equation}
\end{thm}

We denote here the MLE $\widehat{\vartheta}_T^1$, as the optimal input does not depend on $\vartheta$, the MLE reaches efficiency and we deduce its large asymptotic properties.

\begin{thm}\label{sec asymp case 1}
The MLE is uniformly consistent on compacts $\K \subset \mathbb{R}_*^+$, {\it i.e.} for any $\nu >0$,
\begin{equation}\nonumber
\lim_{T\rightarrow \infty} \sup_{\vartheta \in \mathbb{K}}\Pro_\vartheta^{T} \left\{\left|\MLE_T^1-\vartheta \right|> \nu \right\}=0\,,
\end{equation}
uniformly on compacts asymptotically normal: as $T$ tends to $+\infty$,
\begin{equation}\nonumber
\lim_{T\rightarrow\infty}  \sup_{\vartheta \in \mathbb{K}} \left| \Es_\vartheta f\left( \sqrt{T} \left( \hat{\vartheta}_T^1 - \vartheta\right) \right) - \Es f(\xi) \right| =0,\quad \forall f \in {\cal C}_b,
\end{equation}
and $\xi$ is a zero mean Gaussian random variable of variance $\left(\mathcal{I}^1(\vartheta)\right)^{-1}$ (see \eqref{eq:sec fisher 1} for the explicit value)
{\bf which does not depend on $H$}
 and we have the uniform on $\vartheta\in \mathbb{K}$ convergence of the moments: for any $p>0$,
\begin{equation}\nonumber
\lim_{T\rightarrow \infty}  \sup_{\vartheta \in \mathbb{K}} \left| \Es_\vartheta \left|\sqrt{T} \left( \MLE_T^1-\vartheta\right) \right |^p-  \Es \left|\xi\right|^p \right| =0.
\end{equation}
Finally, the MLE is efficient in the sense of \eqref{asymp}.

\end{thm}
\subsection{Case of $k^2 < 2\vartheta$}
In this section, we consider only when $k^2< 2\vartheta$. First of all, we will get the optimal input:
\begin{thm}\label{sec opt case 2}
The asymptotical optimal input in the class of controls $\mathcal{U}_T$ is $u_{opt}^2(t)=\frac{\kappa_H}{\sqrt{2\lambda}}t^{H-\frac{1}{2}}e^{i\omega t}$ where $\omega=\pm \sqrt{\vartheta-\frac{k^2}{2}}$. Moreover,
$$\lim_{T \to + \infty} \frac{{\cal J}_T^2(\vartheta) }{T} =\mathcal{I}^2(\vartheta),$$
where
\begin{equation}\label{eq:sec fisher 2}
\mathcal{I}^2(\vartheta)=\frac{16}{(k^4-4k^2\vartheta)^2}.
\end{equation}
\end{thm}
In this case, the optimal input depends on the parameter $\vartheta$, we can not directly study the properties of MLE, we will use Newton-Raphson method to get the asymptotical properties of MLE which will considered in the Next section.

\section{Preliminary Results}
\subsection{Transformation of The Model}
The explicit representation of the likelihood function can be written thanks to the transformation of observation model proposed in \cite{Klep02b}. In what follows, all random variables and processes are defined on a given stochastic basis $(\Omega,\mathcal{F},(\mathcal{F}_t)_{t \geq 0},\mathbf{P})$ satisfying the usual conditions and processes are $(\mathcal{F}_t)-$ adapted. More over the {\it natural filtration} of a process is understood as the $\mathbf{P}-$ completion of the filtration generated by this process. Let us define:
$$
k_H(t,s)=\kappa_H^{-1}s^{\frac{1}{2}-H}(t-s)^{\frac{1}{2}- H},\,w_H(t)=\frac{1}{2\lambda(2- 2H)}t^{2-2H},
$$
$$
N_t=\int_0^tk_H(t,s)dV_t^H,
$$
where $\kappa_H$ and $\lambda$ are defined in \eqref{sec def_constants}.
Then the process $N=(N_t,t \geq 0)$ is a Gaussian martingale, called in \cite{norrosal} the {\it fundamental martingale}, whose variance function is noting but $w_H$. More over, the natural filtration of the martingale $N$ coincides with the natural filtration of the fBm $V^H$.

Following \cite{Klep02b}, let us introduce a process $Z=(Z_t,0\leq t\leq T)$ the fundamental semi-martingale associated to $Y$, defined as
\begin{equation}\label{sec def semimartingale }
Z_t=\int_0^tk_H(t,s)dY_s.
\end{equation}
Note that $Y$ can be represented as $Y_t=\int_0^tK_H(t,s)dZ_s$, where $K_H(t,s)=H(2H- 1)\int_s^tr^{H- \frac{1}{2}}(r-s)^{H- \frac{3}{2}}dr$ for $0\leq s\leq t$ and therefore the natural filtration of $Y$ and $Z$ coincide. Moreover, we have the following representation:
\begin{equation}\label{sec def semimartingale 2}
dZ_t= \lambda \ell(t)^* \zeta(t)d\langle N\rangle_t+ dN_t,\,Z_0=0,
\end{equation}
where $\zeta(t)$ is the solution of the ordinary differential equation:
\begin{equation}\label{sec tram signal}
\frac{d\zeta(t)}{d\langle N\rangle_t}=\lambda A_0\otimes \mathbf{A}(t)\zeta(t) +b(t)v(t),\,\zeta(0)=0,
\end{equation}
with
\begin{equation}\label{sec matrix}
\ell(t)=\left(
          \begin{array}{c}
            t^{2H-1} \\
            1 \\
            0 \\
            0 \\
          \end{array}
        \right),\, A_0=\left(
                         \begin{array}{cc}
                           0 & 1 \\
                           -\vartheta & -k \\
                         \end{array}
                       \right),\, \mathbf{A}(t)=\left(
                                                  \begin{array}{cc}
                                                    t^{2H-1} & 1 \\
                                                    t^{4H-2} & t^{2H-1} \\
                                                  \end{array}
                                                \right),\, b(t)=\left(
                                                                  \begin{array}{c}
                                                                    0 \\
                                                                    0 \\
                                                                    1 \\
                                                                    t^{2H-1} \\
                                                                  \end{array}
                                                                \right).
\end{equation}
Here, for a control $u(t) $, we define the function $v(t)$ by the following equation:
\begin{equation}\label{sec vtou}
v(t)=\frac{d}{dw_H(t)}\int_0^tk_H(t,s)u(s)ds,
\end{equation}
provided that the fractional derivative exists. Let us define the space of control for $v(t)$ that:
\begin{equation}\label{sec control_spaces}
\mathcal{V}_T=\left\{ v \, \bigg{|} \, \frac{1}{T}\int_0^T|v(t)|^2dw_H(t)\leq 1\right\}.
\end{equation}
Here $|\cdot|$ denote the norm for the complex function. Note that these sets are non empty. Remark that with \eqref{sec vtou} the following relation between control $u$ and its transformation $v$ holds:
\begin{equation}\label{sec eq:utov}
u(t)=\frac{d}{dt}\int_0^tK_H(t,s)v(s)dw_H(s).
\end{equation}
\subsection{Likelihood function and the Fisher information}
The classical Girsanov theorem gives the following equality:
\begin{equation}\label{sec fishe inf}
\mathcal{L}(\vartheta, Z^T)=\exp\left\{\lambda\int_0^T\ell(t)^*\zeta(t)dZ_t- \frac{\lambda^2}{2}\int_0^T\big|\ell(t)^*\zeta(t)\big|^2d\langle N\rangle_t\right\}.
\end{equation}
The fisher information stands for:
$$
\mathcal{I}_T(\vartheta,v)= - \mathbf{E}_{\vartheta}\frac{\partial^2}{\partial \vartheta^2} \ln \mathcal{L}(\vartheta,Z^T),
$$
which is
\begin{equation}\label{sec fishe 5}
\mathcal{I}_T(\vartheta,v)=\lambda^2\int_0^T\left|\ell(t)^*\frac{\partial \zeta(t)}{\partial \vartheta}\right|^2d\langle N\rangle_t.
\end{equation}
\begin{rem}
From the following result we know that for the case $k^2<2\vartheta$, the optimal input depends on the unknown parameter. But in the procedure to find the maximum of the Fisher information we have not consider this situation, that is to say we will only consider the partial derivative of the function $\zeta(t)$ with respect to $\vartheta$ only depends on the function $\varphi(t)$ defined below but not the function $v$.
\end{rem}
\subsection{Proof of Theorem ~\ref{sec opt case 1} and ~\ref{sec opt case 2}}
Let us define
$$
\mathcal{J}(\vartheta)=\sup_{v\in \mathcal{V}_T}\mathcal{I}_T(\vartheta,v).
$$
From \eqref{sec tram signal} we get
\begin{equation}\label{sec signal 2}
\zeta(t)=\varphi(t)\int_0^t\varphi^{-1}(s)b(s)v(s)d\langle N\rangle_s,
\end{equation}
where $\varphi(t)$ is the fundamental matrix satisfying:
\begin{equation}\label{sec varphi}
\frac{d\varphi(t)}{d\langle N\rangle_t}=\lambda A_0\otimes \mathbf{A}(t)\varphi(t),\,\varphi(t)=\mathbf{Id},
\end{equation}
where $\mathbf{Id}$ is the $4\times 4$ identity matrix. Therefore
\begin{equation}\label{sec fishe 2}
\mathcal{I}_T(\vartheta,v)=\int_0^T\int_0^TK_T(s,\sigma)\frac{s^{\frac{1}{2}-H}}{\sqrt{2\lambda}}v(s)
\frac{\sigma^{\frac{1}{2}-H}}{\sqrt{2\lambda}}\bar{v}(\sigma)dsd\sigma,
\end{equation}
where $\bar{v}$ represent the conjugation of the complex function $v$ and
\begin{equation}
K_T(s,\sigma)=\int_{\max(s,\sigma)}^TG(t,s)G(t,\sigma)dt,
\end{equation}
and
\begin{equation}
G(t,\sigma)=\frac{\partial}{\partial\vartheta}\left(\frac{1}{2}t^{\frac{1}{2}-H}\ell(t)^*\varphi(t)\varphi^{-1}(\sigma)b(\sigma)
\sigma^{\frac{1}{2}-H}\right).
\end{equation}
Then
\begin{eqnarray*}
  \mathcal{J}_{T}(\vartheta) &=& T \sup_{\widetilde{v} \in L^2[0,T], \|\widetilde{v}\|\leq 1}\int_0^T\int_0^TK_T(s,\sigma)\widetilde{v}(s)\widetilde{\bar{v}}(\sigma)dsd\sigma ,\\
  &=& T \sup_{\widetilde{v} \in L^2[0,T],\|\widetilde{v}\|\leq 1}(K_T\widetilde{v},\widetilde{\bar{v}}),
\end{eqnarray*}
where $\displaystyle{\tilde{v}(s) = \frac{s^{\frac{1}{2}-H}}{\sqrt{2\lambda}}\frac{v(s)}{\sqrt{T}}}$ and $\|\cdot \|$ stands for the complex norm in $L^2[0,T]$. So in order to prove the Theorem ~\ref{sec opt case 1} and ~\ref{sec opt case 2}, we only need the following two Lemmas.
\begin{lmm}\label{sec cov operator 1}
When $k^2\geq 2\vartheta$,
$$\lim_{T\rightarrow \infty}\sup_{\widetilde{v} \in L^2[0,T],\|\widetilde{v}\|\leq 1}(K_T\widetilde{v},\widetilde{\bar{v}})=\frac{1}{\vartheta^4},$$
with an optimal input $v_{opt}^1(t)=\sqrt{2\lambda}t^{H-\frac{1}{2}}$ belonging to the space of control $\mathcal{V}_T$.
\end{lmm}

\begin{lmm}\label{sec cov operator 2}
When $k^2< 2\vartheta$,
$$\lim_{T\rightarrow \infty}\sup_{\widetilde{v} \in L^2[0,T],\|\widetilde{v}\|\leq 1}(K_T\widetilde{v},\widetilde{\bar{v}})=\frac{16}{(k^4-4k^2\vartheta)^2},$$
with the optimal input $v_{opt}^2(t)=\sqrt{2\lambda}t^{H-\frac{1}{2}}e^{i\omega t}$ where $\omega=\pm \sqrt{\vartheta-\frac{k^2}{2}}$.
\end{lmm}
The proof of Lemma ~\ref{sec cov operator 1} and ~\ref{sec cov operator 2} are based on the Laplace transformation computation and will presented in Section 4.

\subsection{Proof of Theorem ~\ref{sec asymp case 1}}
In order to prove the Theorem~\ref{sec asymp case 1}, we need to check the Ibragimov-Khasminskii Theorem about the asymptotic efficiency for the MLE in \cite{ibragimov}.
\subsubsection{Ibragimov-Khasminskii Theorem}
\begin{thm}\label{asym condition}
Assume that we are given an observable process
$$
d\eta_t=X_T(\vartheta,t)d\langle N\rangle_t+dN_t,\,t\in[0,T],\,\vartheta\in \mathbb{K},
$$

with the following conditions holds:

1.) The Fisher information $\mathcal{I}_T(\vartheta)\rightarrow \infty$ as $T\rightarrow \infty$(or $T\rightarrow 0$) uniformly with respect to $\vartheta\in \mathbb{K}$.

2.) The ratio $\mathcal{I}_T(\vartheta_1)/\mathcal{I}_T(\vartheta_2)$ is uniformly (with respect to T and $\vartheta_i$) bounded.

3.) The function $\vartheta\rightarrow X_T(\vartheta,t)$ is continuously differentiable.

4.) The function $F(h)=\int_0^T|X_T(\vartheta+h\mathcal{I}_T(\vartheta)^{- \frac{1}{2}},t)- X_T(\vartheta,t)|^2d\langle N\rangle_t$ is greater than $C\min(|h|^2,|h|^{\beta})$ for $\vartheta,\,\vartheta+h\mathcal{I}_T(\vartheta)^{-\frac{1}{2}}\in \mathbb{K}$, where $C$ and $\beta$ are positive constants.

Then $\mathcal{I}_T^{1/2}(\vartheta)(\widehat{\vartheta}_T-\vartheta)\Rightarrow \mathcal{N}(0,1)$ as $T\rightarrow \infty$ (or $T\rightarrow 0$), where $\widehat{\vartheta}_T$ is the maximum likelihood estimator for $\vartheta$. Moreover, all moments of $\mathcal{I}_T^{1/2}(\vartheta)(\widehat{\vartheta}_T-\vartheta)$ tend to the corresponding moments of $\mathcal{N}(0,1)$. The convergence is uniform with respect to $\vartheta \in \mathbb{K}$.
\end{thm}
\subsubsection{Taylor's Development Proof}
When $k^2\geq 2\vartheta$, with the optimal input, we can get the new system
\begin{equation}\label{sec changemodel 1}
\left\{\begin{array}{ccrr}
 \frac{d\zeta^1(t)}{d\langle N\rangle_t}&=&\lambda A_0\otimes\mathbf{A}(t)\zeta^1(t)+b(t)v_{opt}^1(t) \,,&\zeta^1(0)=0\,,\\[.2cm]
dZ_t^1 &=& \lambda\ell(t)^*\zeta^1(t)d\langle N\rangle_t+dN_t\,,&Z^1_0=0.
\end{array}\right.
\end{equation}
Let us define the function
\begin{equation}\label{sec def g}
g(\vartheta,t)=t^{\frac{1}{2}}\ell(t)^*\varphi(\vartheta,t)\int_0^t\varphi^{-1}(\vartheta,s)b(s)s^{\frac{1}{2}-H}ds,
\end{equation}
where
$$
\frac{d\varphi(\vartheta,t)}{dt}=\frac{A_0(\vartheta)}{2}\otimes\mathbf{A}_H(t)\varphi(\vartheta,t),
$$
and
$$
\mathbf{A}_H(t)=\left(
                  \begin{array}{cc}
                    1 & t^{1-2H} \\
                    t^{2H-1} & 1 \\
                  \end{array}
                \right).
$$
With Taylor's development with respect to t, we can get that
\begin{equation}\label{Taylor deve}
|g(\vartheta+h,t)-g(\vartheta,t)|=C|h|t^4+o(t^4),
\end{equation}
for every real value $h$, $C$ is a constant which does not depend on $\vartheta$. Here $\frac{o(t^4)}{t^4}=0$ when $t\rightarrow 0$. In our case, the Fisher Information
\begin{eqnarray*}
  \mathcal{I}_T^1(\vartheta,v_{opt}^1 ) &=& \lambda^2\int_0^T\left(\ell(t)^*\frac{\partial \zeta^1(t)}{\partial \vartheta}\right)^2d\langle N\rangle_t \\
   &=& \frac{1}{4}\int_0^T\left|\frac{\partial g(\vartheta,t)}{\partial \vartheta}\right|^2dt,
\end{eqnarray*}
with the condition \eqref{Taylor deve}, we can verify the four conditions in Theorem ~\ref{asym condition}. So that we can get that $\sqrt{T}(\widehat{\vartheta}_T^1-\vartheta)\Rightarrow \mathcal{N}(0, (\mathcal{I}^1(\vartheta))^{-1})$, and moreover, we can get all of the results in the Theorem ~\ref{sec asymp case 1}.

\subsection{Asymptotical Properties of MLE When $k^2<2\vartheta$}
When $k^2<2\vartheta$, with Lemma ~\ref{sec cov operator 2} we know that
$v_{opt}^2(t)=\sqrt{2\lambda}t^{H-\frac{1}{2}}e^{i\omega t}$ where $\omega=\pm \sqrt{\vartheta-\frac{k^2}{2}}$. The optimal input depends on the unknown parameter $\vartheta$. So we can not directly use the Ibragimov-Khasminskii Theorem to find the asymptotic properties of MLE. We follow the general procedure : Divide the observation time interval into two parts, the first one being relatively short. Then find a preliminary estimate $\overline{\vartheta}$ of the unknown parameter $\vartheta$ from the observation in this interval by using an input which does not depend on $\vartheta$. After that, we use $\overline{\vartheta}$, instead of $\vartheta$, to form an approximately optimal input in the second(long) interval. By using this input, we arrive at an asymptotically efficient estimator of $\vartheta$. At the second stage, we can also use the MLE, though this is not an easily-implemented procedure. A simpler method of the Newton-Raphson type can be described as follows.
\subsubsection{Newton-Raphson method}
When $k^2<2\vartheta$, our system is that
\begin{equation}\label{sec changemodel 2}
\left\{\begin{array}{ccrr}
 \frac{d\zeta^2(t)}{d\langle N\rangle_t}&=&\lambda A_0\otimes\mathbf{A}(t)\zeta^2(t)+b(t)v_{opt}^2(t) \,,&\zeta^2(0)=0\,,\\[.2cm]
dZ_t^2 &=& \lambda\ell(t)^*\zeta^2(t)d\langle N\rangle_t+dN_t\,,&Z^2_0=0.
\end{array}\right.
\end{equation}
Here, $Z_t^2$ represents the observable process when we have the $v_{opt}^2(t)$. We know when to find the MLE, we will find the root of the equation
$$
F(Z_t^2,\vartheta)=\int_0^TX_{\vartheta}^{'}(t,\vartheta)dZ_t^2-\int_0^TX(t,\vartheta)X_{\vartheta}^{'}(t,\vartheta)d\langle N\rangle_t=0,
$$
where $X(t,\vartheta)=\lambda \ell(t)^* \zeta^2(t)$ and $X_{\vartheta}^{'}(t,\vartheta)$ is the partial derivative of $X(t,\vartheta)$ with respect to $\vartheta$.

The general Newton iteration method for the solution can described by
$$
\vartheta_{n+1}=\vartheta_n-\frac{F(\vartheta_n)}{F^{'}(\vartheta_n)}.
$$
In fact
$$
F^{'}(Z_t^2,\vartheta)=\int_0^TX_{\vartheta\vartheta}^{''}(t,\vartheta)dZ_t^2-\int_0^T|X_{\vartheta}^{'}(t,\vartheta)|^2d\langle N\rangle_t-\int_0^TX(t,\vartheta)X_{\vartheta\vartheta}^{''}d\langle N\rangle_t.
$$
Or when we develop $dZ_t^2$, we can get that
$$
F^{'}(Z_t^2,\vartheta)=-\int_0^T|X_{\vartheta}^{'}(t,\vartheta)|^2d\langle N\rangle_t+\int_0^TX_{\vartheta\vartheta}^{''}dN_t.
$$
The second term is often negligible compared to the first one. By dropping it and making the first Newton iteration, we get an estimator $\widehat{\vartheta}_T^2$ from an initial estimator $\overline{\vartheta}$ of the parameter $\vartheta$:
$$
\widehat{\vartheta}_T^2=\overline{\vartheta}+\frac{\int_0^TX_{\vartheta}^{'}(t,\overline{\vartheta})dZ_t^2-\int_0^TX(t,\overline{\vartheta})X_{\vartheta}^{'}(t,\overline{\vartheta})d\langle N\rangle_t}{\mathcal{I}(\overline{\vartheta})},
$$
where $\mathcal{I}(\overline{\vartheta})=\mathcal{I}_0^T(\overline{\vartheta})$ is the Fisher information of system \eqref{sec changemodel 2}. The difficulty is that the estimator depends on the observation time $\tau$, the function $X_{\vartheta}^{'}(t,\overline{\vartheta})$ is not non anticipating, the integral $\int_0^{\tau}X_{\vartheta}^{'}(t,\overline{\vartheta})dZ_t^2$, so we can define the estimator as
\begin{equation}\label{estimator1}
\widehat{\vartheta}_T^2=\overline{\vartheta}+\frac{\int_{\tau}^TX_{\vartheta}^{'}(t,\overline{\vartheta})dZ_t^2-\int_{\tau}^TX(t,\overline{\vartheta})X_{\vartheta}^{'}(t,\overline{\vartheta})d\langle N\rangle_t}{\mathcal{I}^T_{\tau}(\overline{\vartheta})},
\end{equation}
or we can write as
\begin{equation}\label{estimator2}
\sqrt{\mathcal{I}_{\tau}^T(\overline{\vartheta})}(\widehat{\vartheta}_T-\vartheta)=\frac{1}{\sqrt{\mathcal{I}_{\tau}^T(\overline{\vartheta})}}
\int_{\tau}^TX_{\vartheta}^{'}(t,\overline{\vartheta})dN_t+R,
\end{equation}
where the remainder
\begin{equation}\label{remainder}
R=\frac{1}{\sqrt{\mathcal{I}_{\tau}^T(\overline{\vartheta})}}\int_{\tau}^TX_{\vartheta}^{'}(t,\overline{\vartheta})\left[
X(t,\vartheta)-X(t,\overline{\vartheta})-X_{\vartheta}^{'}(t,\overline{\vartheta})(\vartheta-\overline{\vartheta})\right]d\langle N\rangle_t,
\end{equation}
with Taylor formula, we can write it as
\begin{equation}\label{remainder2}
R=\frac{1}{\sqrt{\mathcal{I}_{\tau}^T(\overline{\vartheta})}}\int_{\tau}^TX_{\vartheta}^{'}(t,\overline{\vartheta})X_{\vartheta\vartheta}^{''}(t,\overline{\overline{\vartheta}})d\langle N\rangle_t(\vartheta-\overline{\vartheta})^2,
\end{equation}
where $\overline{\overline{\vartheta}}$ is a point between $\vartheta$ and $\overline{\vartheta}$. In view of equation \eqref{estimator2}, we know that, if there is no the remainder $R$, we can get the asymptotic efficiency of the estimator $\widehat{\vartheta}_T^2$ since $\mathcal{I}_{\tau}^T(\overline{\vartheta})$ is asymptotically equivalent to $\mathcal{I}_T^2(\vartheta, v_{opt}^2(t))$ which is the Fisher information of the system \eqref{sec changemodel 2}. So we need the remainder is small. To study the remainder, we need to study the estimator $\overline{\vartheta}$.
\subsubsection{Small Interval Estimator}
We will observe the small interval $[0, \tau]$. Let us define a function $\overline{v}(t)=\rho\sqrt{2\lambda}t^{H-\frac{1}{2}}$ where $\rho$ is a constant depending on $\tau$. Assume that we are given a linear system
\begin{equation}\label{sec changemodel 3}
\left\{\begin{array}{ccrr}
 \frac{d\zeta^3(t)}{d\langle N\rangle_t}&=&\lambda A_0\otimes\mathbf{A}(t)\zeta^3(t)+b(t)\overline{v}(t) \,,&\zeta^3(0)=0\,,\\[.2cm]
dZ_t^3 &=& \lambda\ell(t)^*\zeta^3(t)d\langle N\rangle_t+dN_t\,,&Z^3_0=0,
\end{array}\right.
\end{equation}
where $Z_t^3$ is the observable process and we only observe the interval $[0, \tau]$ and get the estimator which we define as  $\overline{\vartheta}$. We have the following Lemma:
\begin{lmm}\label{small efficient}
when give the system \eqref{sec changemodel 3}, the MLE  $\overline{\vartheta}$ for the parameter $\vartheta\in \mathbb{K}$ is asymptotically  efficient provided that when $\tau\rightarrow 0$, $\tau^9\rho^2\rightarrow \infty$.
\end{lmm}
\paragraph{The proof} Follows from the Ibragimov-Khasminskii Theorem, we know that the Fisher information is
$$
\mathcal{I}_{\tau}(\vartheta)=\frac{1}{4}\int_0^{\tau}\left|\rho\frac{\partial}{\partial \vartheta}g(\vartheta, t)\right|^2dt
$$
where $g(\vartheta, t)$ is defined in \eqref{sec def g}. With Taylor's development of $g(\vartheta,t)$ and the condition $\tau^9\rho^2\rightarrow \infty$, we can verify the four conditions of Ibragimov-Khasminskii Theorem when $\tau\rightarrow 0$. So this Lemma follows.

In fact, we can have a more advanced result:
\begin{corr}\label{corr1}
If we choose an arbitrary small interval $[0, \tau]$, $\tau=o(1)$, as $T\rightarrow \infty$, and
$$
\int_0^{\tau}\left|\overline{v}(t)\right|^2d\langle N\rangle_t=o(T),
$$
we can obtain the estimator with the precision of order $\frac{1}{\sqrt{T}}$. More precisely, if $f(T)=o(\sqrt{T})$, we can find the estimator $\overline{\vartheta}$ such that $\overline{\vartheta}-\vartheta=O\left(\frac{1}{f(T)}\right)$.
\end{corr}
\subsubsection{Long Time Estimation}
Now, we will return to the system \eqref{sec changemodel 2}, but using the estimator $\overline{\vartheta}$. we define a new function $\underline{v}(t)$ to replace $v_{opt}^2(\overline{\vartheta}, t)$, that is
$$
\underline{v}(t)=\sqrt{2\lambda}t^{H-\frac{1}{2}}e^{\pm i \sqrt{\overline{\vartheta}-\frac{k^2}{2}}},
$$
and get the MLE defined in \eqref{estimator2}. We have the following Theorem:
\begin{thm}
When given the system \eqref{sec changemodel 2}, we have an asymptotically efficient two-stage estimator of parameter $\vartheta$. The first is given in Corollary ~\ref{corr1} (where, e.g. $\tau=T^{-\varepsilon}$, $\rho=\sqrt{T}$, $\varepsilon$ is small positive). The second stage is $\widehat{\vartheta}_T^2$ defined in \eqref{estimator2}.
\end{thm}

\section{Proof of Lemma ~\ref{sec cov operator 1} and Lemma ~\ref{sec cov operator 2}}
In this section, we will prove Lemma ~\ref{sec cov operator 1} and Lemma ~\ref{sec cov operator 2}. First of all, we will use Laplace transform to get the upper bound of the operator $K_T$. Then we will get the lower bound using the special value of the function $v(t)$.
\begin{rem}
Even the input $v$ can be a complex function, but the operator $K_T(s,\sigma)$ is still a real symmetric operator, so the method to find the upper bound in \cite{BKP10} can still be used in our situation.
\end{rem}
\subsection{Laplace Transform Proof of Upper Bound}
Let us introduce the pair process $\xi=((\xi_1,\xi_2),0\leq t\leq T)$ with

\begin{equation}\label{sec process xi1}
\xi_t^1=\left(\int_t^T\sigma^{\frac{1}{2}-H}\ell(\sigma)^*\varphi(\sigma)*dW_{\sigma}\right)\varphi^{-1}(t),
\end{equation}
and
\begin{equation}\label{sec process xi2}
\xi_t^2=\frac{\partial}{\partial \vartheta}\xi_t^1,
\end{equation}
where $W$ is a Winer process and $*dW_{\sigma}$ denotes the It\^{o} backward integral (see \cite{Roz}). It is worth emphasizing that
$$
K_T(s,\sigma)=\frac{1}{4}\mathbf{E}\left(\xi_s^2b(s)s^{\frac{1}{2}-H}\xi^2_{\sigma}b(\sigma)\sigma^{\frac{1}{2}-H}\right)=\mathbf{E}(\mathcal{X}_{\sigma}\mathcal{X}_s),
$$
where $\mathcal{X}$ is the centered Gaussian process defined by :
$$
\mathcal{X}_t=\frac{1}{2}\xi_t^2b(t)t^{\frac{1}{2}-H}.
$$
The process $\xi$ also satisfies the following dynamic:
$$
-d\xi_t=\xi_t\mathcal{A}(t)d\langle N\rangle_t+\mathcal{L}(t)*dM_t,\,\xi_T=0,
$$
with $M=(M_t,\,t\geq 0)$ a martingale of the same variance function as $N=(N_t,\,t\geq 0)$,

$\mathcal{A}(t)=\left(
                  \begin{array}{cccc}
                    0 & 1 & 0 & 0 \\
                    -\vartheta & - k & -1 & 0 \\
                    0 & 0 & 0 & 1 \\
                    0 & 0 & -\vartheta & -k \\
                  \end{array}
                \right)
\otimes \lambda \mathbf{A}(t)$ and $\mathcal{L}(t)=\sqrt{2\lambda}\left(
                                                                    \begin{array}{cc}
                                                                      \ell(t)^* & 0 \\
                                                                    \end{array}
                                                                  \right).
$

In fact, the covariance operator $K_T$ is a symmetrical compact operator, we should estimate the spectral gap(the first eigenvalue $\nu_1(T)$). This estimation is based on the Laplace transform computation.
Let us compute, for sufficiently small negative $a<0$ the Laplace transform of $\int_0^T\mathcal{X}_t^2dt$:
\begin{eqnarray*}
  L_T(a) &=& \mathbf{E}_{\vartheta}\exp\left\{-a\int_0^T\mathcal{X}_t^2dt\right\} \\
   &=& \mathbf{E}_{\vartheta}\exp\left\{- a\int_0^T\left[\frac{1}{2}\left(\frac{\partial}{\partial \vartheta}\xi_t^1\right)b(t)t^{\frac{1}{2}-H}\right]^2dt\right\}.
\end{eqnarray*}
On the one hand, for $a>- \frac{1}{\nu_1(T)}$, since $\mathcal{X}$ is a centered Gaussian process with covariance operator $K_T$, using Mercer's theorem and Parseval's inequality, $L_T(a)$ can be represented as:
\begin{equation}\label{sec laplace}
L_T(a)=\prod_{i\geq 1}(1+2a\nu_i(T))^{-\frac{1}{2}},
\end{equation}
where $\nu_i(T)$, $i\geq 1$ is the sequence of positive eigenvalues of the covariance operator. On the other hand,
\begin{eqnarray*}
  L_T(a) &=& \mathbf{E}_{\vartheta}\exp\left\{-\frac{a\lambda}{2}\int_0^T\xi_t\mathcal{M}\xi_t^*d\langle N\rangle_t  \right\} \\
   &=& \exp\left\{\frac{1}{2}\tr(\mathcal{H}(t)\mathcal{L}(t)^*\mathcal{L}(t))d\langle N\rangle_t\right\},
\end{eqnarray*}
where
$$
\mathcal{M}(t)=\left(
                 \begin{array}{cc}
                   0 & 0 \\
                   0 & b(t)b(t)^* \\
                 \end{array}
               \right),
$$
and $\mathcal{H}(t),\,t\geq 0$ is the solution of Ricatti differential equation:
\begin{equation}\label{sec Rica eq}
\frac{d\mathcal{H}(t)}{d\langle N\rangle_t}=\mathcal{H}(t)\mathcal{A}(t)^*+\mathcal{A}(t)\mathcal{H}(t)+\mathcal{H}(t)\mathcal{L}(t)^*\mathcal{L}(t)\mathcal{H}(t)- a \lambda \mathcal{M}(t),
\end{equation}
with initial value $\mathcal{H}(0)=0$, provided that the solution of equation $\eqref{sec Rica eq}$ exists for any $0\leq t\leq T$.
It is well known that if $\det \Psi_1(t)>0$, for any $t\in [0,T]$, then the solution $\mathcal{H}$ of equation \eqref{sec Rica eq} can be written as $\mathcal{H}(t)=\Psi_1^{-1}(t)\Psi_2(t)$, where the pair of $8\times 8$ matrices $(\Psi_1,\Psi_2)$ satisfies the system of linear differential equation:
\begin{equation}\label{Psi1}
\frac{d\Psi_1(t)}{d\langle N\rangle_t}=-\Psi_1(t)\mathcal{A}(t)- \Psi_2(t)\mathcal{L}(t)^*\mathcal{L}(t),\,\Psi_1(0)=\mathbf{Id}_{8\times 8},
\end{equation}
$$
\frac{d\Psi_2(t)}{d\langle N\rangle_t}-a\lambda\Psi_1(t)\mathcal{M}(t)+\Psi_2(t)\mathcal{A}(t)^*,\,\Psi_2(0)=\mathbf{0}.
$$

Moreover, under the condition $\det \Psi_1(t)>0$, for any $t\in[0,T]$, the following equality holds:
\begin{eqnarray*}
  L_T(a) &=& \exp\left\{-\frac{1}{2}\int_0^Ttrace\mathcal{A}(t)d\langle N\rangle_t\right\}(\det\Psi_1(T))^{-\frac{1}{2}} \\
   &=& \exp \{kT\}(\det\Psi_1(T))^{-\frac{1}{2}},
\end{eqnarray*}
or equivalently using \eqref{sec laplace},
\begin{equation}\label{sec laplace 1}
\prod_{i\geq 1}(1+2a\nu_i(T))=\exp\{-2kT\}(\det\Psi_1(T)).
\end{equation}

Let us note here that the solution of linear system \eqref{Psi1} exist for any $t>0$ and for any $a\in \mathbb{C}$. For $a=0$, $\det \Psi_1(t)=\exp\{2kt\}>0$. Due to the continuity property of the solutions of linear differential equations with respect to a parameter, for all $T>0$, there exists $a(T)<0$ such that
$$
\inf_{t\in[0,T]}\det \Psi_1(t)>0.
$$
Therefore, equality \eqref{sec laplace 1} holds in an open set in $\mathbb{C}$, containing 0. Compactness of the covariance operator implies due to the Weierstrass theorem, the analytic property of $\underset{i\geq 1}{\prod}(1+2a\nu_i(T))$ with respect to $a$. Hence, equality \eqref{sec laplace 1} holds for any $a\in \mathbb{C}$.

Now, we rewrite the system of $(\Psi_1,\Psi_2)$ that
\begin{equation}\label{sec psi to}
\frac{d(\Psi_1(t),\Psi_2(t)\mathbf{J})}{d\langle N\rangle_t}=(\Psi_1(t),\Psi_2(t)\mathbf{J})\cdot (\Upsilon \otimes \lambda \mathbf{A}(t))
\end{equation}
where $\mathbf{J}=\left(
                    \begin{array}{cccc}
                      J & J & J & J \\
                      J & J & J & J \\
                      J & J & J & J \\
                      J & J & J & J \\
                    \end{array}
                  \right)
$ and $J=\left(
         \begin{array}{cc}
           0 & 1 \\
           1 & 0 \\
         \end{array}
       \right)
$  and
$$
\Upsilon=\left(
  \begin{array}{cccccccc}
    0 & -1 & 0 & 0 & 0 & 0 & 0 & 0 \\
    \vartheta & k & 1 & 0 & 0 & 0 & 0 & 0 \\
    0 & 0 & 0 & -1 & 0 & 0 & 0 & 0 \\
    0 & 0 & \vartheta & k & 0 & 0 & 0 & -a \\
    -2 & 0 & 0 & 0 & 0 & - \vartheta & 0 & 0 \\
    0 & 0 & 0 & 0 & 1 & -k & 0 & 0 \\
    0 & 0 & 0 & 0 & 0 & -1 & 0 & -\vartheta \\
    0 & 0 & 0 & 0 & 0 & 0 & 1 & -k \\
  \end{array}
\right).
$$

The eigenfunction of $\Upsilon$ is that
$$
(y^2-ky+\vartheta)^2(y^2+ky+\vartheta)^2+2a=0.
$$
\subsubsection{The Case Of $k^2\geq 2\vartheta$}
In this case, when $-\frac{\vartheta^4}{2}<a<0$, let $(y_i)_{i=1,\ldots 8}$ be the eigenvalues of the matrix $\Upsilon$, it can be checked that
$$
\det\Psi_1(T)=\exp\left((y_1+y_3+y_5+y_7)T\right)(C+O(\frac{1}{T})).
$$
where C is a constant and there are 3 cases with different $y_i$.

\paragraph{(1) $k^2\geq 4\vartheta$ , there are 8 real eigenvalues} we get that
$$
y_1=\sqrt{\frac{k^2-2\vartheta+\sqrt{k^4-4k^2\vartheta+4\sqrt{-2a}}}{2}},
$$

$$
y_3=\sqrt{\frac{k^2-2\vartheta-\sqrt{k^4-4k^2\vartheta+4\sqrt{-2a}}}{2}},
$$

$$
y_5=\sqrt{\frac{k^2-2\vartheta+\sqrt{k^4-4k^2\vartheta-4\sqrt{-2a}}}{2}},
$$

$$
y_7=\sqrt{\frac{k^2-2\vartheta-\sqrt{k^4-4k^2\vartheta-4\sqrt{-2a}}}{2}}.
$$

\paragraph{(2) $k^2\geq 4\vartheta$ or $2\vartheta\leq k^2 \leq 4\vartheta$, there are 4 real eigenvalues and 4 complex eigenvalues}
$$
y_1=\sqrt{\frac{k^2-2\vartheta+\sqrt{k^4-4k^2\vartheta+4\sqrt{-2a}}}{2}},
$$

$$
y_3=\sqrt{\frac{k^2-2\vartheta-\sqrt{k^4-4k^2\vartheta+4\sqrt{-2a}}}{2}},
$$

$y_5+y_7=2\sqrt{m}$ where $m$ and $n$ are the solutions of the equation $m^2-n^2=k^2-2\vartheta$ and $2mn=\sqrt{4\sqrt{-2a}-k^4+4k^2\vartheta}$.

\paragraph{(3) $2\vartheta\leq k^2 \leq 4\vartheta$, there are 8 complex eigenvalues}

$y_1+y_3=2\sqrt{p}$ where $p$ and $q $ are the solutions of the equation $p^2-q^2=k^2-2\vartheta$ and $2pq=\sqrt{-4\sqrt{-2a}-k^4+4k^2\vartheta}$. $y_5+y_7=2\sqrt{m}$ where $m$ and $n$ are the solutions of the equation $m^2-n^2=k^2-2\vartheta$ and $2mn=\sqrt{4\sqrt{-2a}-k^4+4k^2\vartheta}$.

Therefore, due to the equality \eqref{sec laplace 1}, we have that when $k^2\geq 2\vartheta$, $\underset{i\geq 1}{\prod}(1+2a\nu_i(T))>0$ for any $a>-\frac{\vartheta^4}{2}$. It means that
$$
\nu_1(T)\leq \frac{1}{\vartheta^4}.
$$

\subsubsection{The Case Of $k^2 <2\vartheta$}
Now let us consider $k^2<2\vartheta$, when $-\frac{(k^4-4k^2\vartheta)^2}{32}<a<0$, there are 8 complex eigenvalues and  it can be check that
$$
\det\Psi_1(T)=\exp\left((y_1+y_3+y_5+y_7)T\right)(C+O(\frac{1}{T})),
$$
where  $y_1+y_3=2\sqrt{p}$ , $p$ and $q $ are the solutions of the equation $p^2-q^2=k^2-2\vartheta$ and $2pq=\sqrt{-4\sqrt{-2a}-k^4+4k^2\vartheta}$.  $y_5+y_7=2\sqrt{m}$ , $m$ and $n$ are the solutions of the equation $m^2-n^2=k^2-2\vartheta$ and $2mn=\sqrt{4\sqrt{-2a}-k^4+4k^2\vartheta}$.

Therefore, with  the equality \eqref{sec laplace 1}, $\underset{i\geq 1}{\prod}(1+2a\nu_i(T))>0$ for any $a>-\frac{(k^4-4k^2\vartheta)^2}{32}$ which means that
$$
\nu_1(T)\leq \frac{16}{(k^4-4k^2\vartheta)^2}.
$$

\subsection{Lower Bound of The Operator}
For the lower bound we only need to calculate the
\begin{equation}
\lim_{T\rightarrow \infty}\frac{\lambda^2}{T}\int_0^T \left(\frac{\partial \zeta^o(t)}{\partial \vartheta}\right)^*\ell(t)\ell(t)^*\left(\frac{\partial \zeta^o(t)}{\partial \vartheta}\right)
\end{equation}
where
$$
\frac{d\zeta^o(t)}{d\langle N\rangle_t}=\lambda A_0\otimes \mathbf{A}(t)\zeta^o(t) +b(t)v_{opt}(t),\,\zeta(0)=0.
$$
The computation will be the same as in \cite{BC13} and in the section 2.5 of \cite{BKP10} we have a important result that for $t$ and $s$ large enough:
$$
g(t,s) \sim 2e^{-\vartheta(t-s)}+\frac{(2H-1)^4}{2\vartheta^2ts},
$$
where
$$
g(t,s)=t^{1/2-H}\left(
                  \begin{array}{c}
                    t^{2H-1} \\
                    1 \\
                  \end{array}
                \right)^* \alpha(t)\alpha^{-1}  \left(
                              \begin{array}{c}
                                1 \\
                                s^{2H-1} \\
                              \end{array}
                            \right)s^{1/2-H},
$$
and the deterministic equation $\alpha(t)$ is defined in the equation (17) in the article \cite{BKP10}. When we compute the limit result, the part of $\frac{(2H-1)^4}{2\vartheta^2ts}$ will be $0$. So the lower bound will be the same in the model driven by the standard Brownian motion which is
\begin{equation}\label{eq:Brownian motion}
\left\{\begin{array}{ccrr}
 \frac{d^2x}{dt^2}+k\frac{dx}{dt}+\vartheta x&=&u(t) \,,&x(0)=0\,,\\[.2cm]
dY_t &=& x(t)dt+dW_t\,,&Y_0=0,
\end{array}\right.
\end{equation}
where $W_t$ is a standard Brownian motion. We can get the Fisher information of this system
\begin{equation}\label{sec fishe 3}
\mathcal{I}_T(\vartheta, u)=\int_0^T\left(\frac{\partial x(t)}{\partial \vartheta}\right)^2dt.
\end{equation}
Or we can write as
\begin{equation}\label{sec fishe 4}
\mathcal{I}_T(\vartheta, u)=\int_0^T\left[\left(
                                            \begin{array}{cc}
                                              1 & 0 \\
                                            \end{array}
                                          \right)\frac{\partial X(t)}{\partial \vartheta}
        \right]^2dt,
\end{equation}
where $X(t)$ is the solution of the equation
$$
dX(t)=A_0X(t)dt+\left(
                  \begin{array}{c}
                    0 \\
                    1 \\
                  \end{array}
                \right)u(t)dt,
$$
$A_0$ is defined in \eqref{sec matrix}. The result in \cite{Oss} tells us
$$
\frac{\mathcal{I}_T(\vartheta, v_{opt}^1(t))}{T}=\frac{1}{\vartheta^4},
$$
and
$$
\frac{\mathcal{I}_T(\vartheta, v_{opt}^2(t))}{T}=\frac{16}{(k^4- 4k^2\vartheta)^2}.
$$
Which achieves the proof.







\end{document}